\documentclass{ws-ijbc}
\usepackage{graphicx}
\usepackage{epstopdf}

\usepackage{amsfonts}
\usepackage{amscd}
\newcommand\Z{\mathbb Z}
\newcommand\R{\mathbb R}
\newcommand\T{\mathbb T}

\newcommand\ee{{\rm e}}

\newcommand\eps\varepsilon

\newcommand\J{\mathcal J}

\newcommand\Lc{\mathcal L}
\newcommand\M{\mathcal M}

\newcommand\p[1]{\left(#1\right)}
\newcommand\pq[1]{\left[#1\right]}
\newcommand\pp[1]{\left\{#1\right\}}
\newcommand\scprod[2]{\left\langle#1,#2\right\rangle}
\newcommand\abs[1]{\left|#1\right|}

\newcommand\norm[1]{\left\|#1\right\|}

\newcommand\tl{\tilde}
\newcommand\wh[1]{\widehat{#1}}
\newcommand\Ord{{\mathcal O}}
\newcommand\rint{\mathop{\rm rint}}
\newcommand\hardsum[2]
  {\!\!\!\!\ds\sum_{\begin{array}{c}\null^{#1}\\[-3pt]\null^{#2}\end{array}}\!\!\!}

\newcommand\mmatrix[4]{\p{\begin{array}{cc}#1&#2\\[3pt]#3&#4\end{array}}}
\newcommand\symmatrix[3]{\mmatrix{#1}{#2}{#2}{#3}}

\newcommand\bremarks{\bigskip\noindent{\bf Remarks.}\ \medskip\bnm}
\newcommand\eremarks{\enm\bigskip\medskip}
\newcommand\qed{\ \ \null\nolinebreak\hfill$\frame{\large\phantom a}$}
\newcommand\paragr[1]{\ \\[12pt]\noindent\textbf{#1.}\quad}

\newcommand\beq{\begin{equation}}
\newcommand\eeq{\end{equation}}
\newcommand\bea{\begin{eqnarray}}
\newcommand\eea{\end{eqnarray}}
\newcommand\bean{\begin{eqnarray*}}
\newcommand\eean{\end{eqnarray*}}
\newcommand\btm{\vspace{-\baselineskip}\begin{itemize}}
\newcommand\etm{\end{itemize}\vspace{-\baselineskip}}
\newcommand\btmm{\begin{itemize}}
\newcommand\etmm{\end{itemize}}
\newcommand\bnm{\vspace{-\baselineskip}\begin{enumerate}}
\newcommand\enm{\end{enumerate}\vspace{-\baselineskip}}

\begin{document}

\title{Exponentially small lower bounds for the splitting of
  separatrices to whiskered tori with frequencies of constant type
  \ \footnote{This work has been partially supported by the Spanish
      MINECO-FEDER Grants MTM2009-06973, MTM2012-31714
      and the Catalan Grant 2009SGR859.
      The author~MG has also been supported by the
      DFG~Collaborative Research Center TRR~109
      ``Discretization in Geometry and Dynamics''.}}
\author{\sc
    Amadeu Delshams$\,^1$,
  \ Marina Gonchenko$\,^2$,
  \\\sc
  \ Pere Guti\'errez$\,^1$
\\[12pt]
  {\small
  $^1\;$\parbox[t]{6.2cm}{
    Dep. de Matem\`atica Aplicada I\\
    Universitat Polit\`ecnica de Catalunya\\
    Av. Diagonal 647, 08028 Barcelona\\
    {\footnotesize
      \texttt{amadeu.delshams@upc.edu}\\
      \texttt{pere.gutierrez@upc.edu}}}
  \quad
  $^2\;$\parbox[t]{6.2cm}{
    Institut f\"ur Mathematik, MA 7-2\\
    Technische Univesit\"at Berlin\\
    Stra{\ss}e des 17. Juni 136, 10623 Berlin\\
    {\footnotesize
      \texttt{gonchenk@math.tu-berlin.de}}}
  }}
\maketitle

\begin{abstract}
We study the splitting of invariant manifolds of whiskered tori with two
frequencies in nearly-integrable Hamiltonian systems, such that the hyperbolic
part is given by a pendulum. We consider a 2-dimensional torus with a fast
frequency vector $\omega/\sqrt\eps$, with $\omega=(1,\Omega)$ where $\Omega$ is
an irrational number of constant type, i.e.~a number whose continued fraction
has bounded entries. Applying the Poincar\'e--Melnikov method, we find
exponentially small lower bounds for the maximal splitting distance between the
stable and unstable invariant manifolds associated to the invariant torus, and
we show that these bounds depend strongly on the arithmetic properties of the
frequencies.
\end{abstract}


\section{Introduction}

The aim of this paper is to introduce a methodology for measuring the
exponentially small splitting of separatrices in a perturbed Hamiltonian
system, associated to a 2-dimensional whiskered torus
(invariant hyperbolic torus), with \emph{fast frequencies}:
\begin{equation}
  \omega_\varepsilon = \frac\omega{\sqrt{\varepsilon}}\,,
  \quad \omega=(1,\Omega),
  \quad\eps>0.
\label{eq:omega_eps}
\end{equation}
This phenomenon requires a careful study due
to the \emph{singular} character of the problem with respect to the
perturbation parameter $\eps$ and, on the other hand,
to the presence of \emph{small divisors}
associated to the frequencies of the torus.
For this reason, the first
results on asymptotic estimates for the splitting are recent,
and have been obtained assuming concrete frequency ratios $\Omega$,
such as the golden mean, and other quadratic irrational numbers.
In this paper, we show that a partial generalization of such techniques
allows us to obtain \emph{lower bounds} for the splitting,
for \emph{all} sufficiently small values of $\eps$,
assuming that the frequency ratio is an
irrational \emph{number of constant type}
(also called a badly approximable number),
i.e.~a number whose continued fraction has bounded entries.
In this way, we establish the existence of splitting for a
much wider (uncountable) class of frequency ratios.

As the unperturbed system, we consider an integrable Hamiltonian $H_0$
with $3$ degrees of freedom having $2$-dimensional whiskered tori
with coincident stable and unstable whiskers (invariant manifolds).
In general, for a perturbed Hamiltonian
\beq\label{eq:Hpert}
  H=H_0+\mu H_1
\eeq
where $\mu$ is small, the stable and unstable whiskers
of a given torus do not coincide anymore,
giving rise to the phenomenon called \emph{splitting of separatrices},
discovered by Poincar\'e \cite{Poincare90}.
In order to give a measure for the splitting,
one often describes it by a periodic vector function
$\M(\theta)$, $\theta\in\T^2$,
usually called \emph{splitting function},
giving the distance between the invariant manifolds
in the complementary directions, on a transverse section.
The most popular tool to measure the splitting is
the \emph{Poincar\'e--Melnikov method},
introduced in \cite{Poincare90}
and rediscovered later \cite{Melnikov63,Arnold64}.
This method provides a first order approximation
\beq\label{eq:Melniapprox}
  \M(\theta)=\mu M(\theta)+\Ord(\mu^2),
\eeq
where $M(\theta)$ is called the \emph{Melnikov function}
and is defined by an integral.
In fact, it was established \cite{Eliasson94,DelshamsG00} that both
vector functions are the gradients of
scalar functions: the \emph{splitting potential} and the
\emph{Melnikov potential}, denoted $\Lc(\theta)$ and $L(\theta)$ respectively.
This result implies the existence of \emph{homoclinic orbits}
(i.e.~intersections between the stable and unstable whiskers)
in the perturbed system.

We consider in~(\ref{eq:Hpert}) a whiskered torus with fast frequencies
as in~(\ref{eq:omega_eps}),
assuming a relation between the parameters $\eps$ and $\mu$
of the form $\mu=\eps^p$ for some $p>0$, having in this way a singular problem.
The interest for such a setting
lies in its relation to the normal form of a nearly-integrable Hamiltonian,
with $\eps$ as the perturbation parameter,
in the vicinity of a simple resonance \cite{Niederman00,DelshamsG01}.

In such a singular problem, one can show that the splitting
is \emph{exponentially small} with respect to $\eps$.
The first results on exponentially small splitting concerned the case
of one and a half degrees of freedom, i.e.~for 1~frequency,
providing upper bounds \cite{Neishtadt84}.
The problem of establishing \emph{asymptotic estimates},
or at least \emph{lower bounds},
for the exponentially small splitting, is more difficult,
due to the fact that the Melnikov function is exponentially small in $\eps$
and the error of the method could overcome the main term
in~(\ref{eq:Melniapprox}).
The first result justifying the Poincar\'e--Melnikov method
and, hence, providing an asymptotic estimate for the
exponentially small splitting was obtained in \cite{Lazutkin03}
for the Chirikov standard map. Later, this was extended
to the case of a Hamiltonian with one and a half degrees of freedom
\cite{DelshamsS92,DelshamsS97,Gelfreich97}
or an area-preserving map \cite{DelshamsR98}.
In the quoted papers specific perturbations $H_1$ were assumed,
but a more general (meromorphic) perturbation
was recently considered in \cite{GuardiaS12}.
It is worth remarking that, in some cases,
the Poincar\'e--Melnikov method does not predict correctly
the size of the splitting, as shown for instance in \cite{BaldomaFGS12}.

For~2 or more frequencies, it turns out that
\emph{small divisors} appear in the splitting function and,
as first noticed in \cite{Lochak92},
the arithmetic properties of the frequency vector $\omega$ play an
important r\^ole. This was established in \cite{Simo94},
and rigorously proved in \cite{DelshamsGJS97}
for the quasi-periodically forced pendulum.
A different technique was used in \cite{LochakMS03}
(see also \cite{RudnevW00}),
namely the parametrization of the whiskers as solutions of
Hamilton--Jacobi equation, to obtain exponential small estimates
of the splitting, and the existence of transverse homoclinic orbits
for some intervals of the perturbation parameter~$\varepsilon$.
Moreover, it was shown in \cite{DelshamsG04} the continuation of
the exponentially small estimates and the
transversality of the splitting, for all sufficiently small values of $\eps$,
under a certain condition on the phases of the perturbation.
Otherwise, homoclinic bifurcations can occur,
studied in \cite{SimoV01} for the Arnold's example.
The quoted papers considered the case of 2~frequencies,
and assuming in most cases that the frequency ratio
is the famous \emph{golden mean}: $\Omega=(\sqrt5-1)/2$.
A~generalization to some other 2-dimensional quadratic frequencies as well to
3-dimensional cubic frequencies was studied
in \cite{DelshamsG03,DelshamsGG13},
and the case of 2-dimensional frequencies of constant type had also
been considered in \cite{RudnevW98}.
For a more complete background and references concerning exponentially small
splitting, see for instance {\cite{LochakMS03,DelshamsGS04}}.

In the present paper we pay attention
to non-quadratic frequency vectors (\ref{eq:omega_eps})
of constant type. It turns out
that the methodology developed in \cite{DelshamsG03,DelshamsGG13}
can be partially applied and that,
using arithmetic properties of numbers of constant type,
we obtain exponentially small \emph{lower bounds} for the
\emph{maximal splitting distance}
(and, consequently, show the existence of splitting).

It is well-known that the property of being a number of constant type
is equivalent to satisfying a \emph{Diophantine condition}
with the minimal exponent (see, for instance, \cite{Lang95}):
there exists $\nu>0$ such that
\begin{equation}\label{eq:DiophcondC}
|q \Omega - p|\ge\frac{\nu}{q},
  \qquad
  \forall p, q\in \Z, \ q\ge1.
\end{equation}
This condition can also be expressed in terms of the vector $k=(-p,q)$:
for some $\gamma>0$,
\begin{equation}
|\langle k, \omega\rangle| \geq \frac{\gamma}{|k|}, \;\;\;
\forall k\in \mathbb{Z}^2\setminus\{0\}.
\label{eq:DiophCond}
\end{equation}
One of the goals of this paper is to show, for the above frequencies,
that we can detect the integer vectors
$k\in\Z^2\setminus\pp0$ that fit better the inequality
in~(\ref{eq:DiophCond}). We show that the ``least'' small divisors
(relatively to the size of $\abs k$) are related to
\emph{principal convergents} $p_n/q_n$ of the continued fraction of $\Omega$,
and we call such vectors $k=v(n)=(-p_n,q_n)$
the \emph{resonant convergents} of $\omega$.
This allows us to detect the dominant harmonic
in the splitting function $\M(\theta)$,
for each small enough value of the perturbation parameter $\eps$.

In the main result of this paper (see Theorem~\ref{thm:main}),
we establish exponentially small lower bounds
for the \emph{maximal splitting distance}, valid in
the case of frequencies of constant type,
giving in this way a partial generalization
of the asymptotic estimates obtained in \cite{DelshamsGG13}
for quadratic frequencies. We point out that the set of
irrational numbers of constant type is uncountable,
in contrast to quadratic irrational numbers that form a countable set.
We also stress that, for some purposes, it is not necessary
to establish the transversality of the splitting, and can be enough
to provide lower bounds of the maximal splitting distance.
Indeed, such lower bounds imply the existence of splitting between the
invariant manifolds, which provides a strong indication
of the non-integrability of the system near the given torus,
and opens the door to the application
of topological methods \cite{GideaL06} for the study of
Arnold diffusion in such systems.

\subsection{Setup and main result}

In order to formulate our main result,
let us describe the Hamiltonian considered,
which is analogous to the one considered
in \cite{DelshamsGS04} and other related works.
In symplectic coordinates
$(x,y,\varphi,I)\in\T\times\R\times\T^2\times\R^2$,
we consider a perturbed Hamiltonian~(\ref{eq:Hpert}),
with $H_0$, $H_1$ of the form
\bea
  \label{eq:ham1}
  &&H_0 (x, y, I) =
  \scprod{\omega_\eps}I+\frac12\scprod{\Lambda I}I+\frac{y^2}2+\cos x-1,
\\
  \label{eq:ham2}
  &&H_1 (x, \varphi)= h(x)\,f(\varphi),
\\
  \label{eq:ham3}
  &&h(x) = \cos x,
  \qquad
  f(\varphi)=\hardsum{k\in \mathbb{Z}^2}{k_2\ge0}
  \ee^{-\rho |k|} \cos(\langle k, \varphi \rangle - \sigma_k),
\eea
where the restriction in the sum is introduced in order to avoid repetitions.
This Hamiltonian is a generalization of the Arnold example
(introduced in \cite{Arnold64} to illustrate the transition
chain mechanism in Arnold diffusion).
It provides a particular model for the behavior
of a nearly-integrable Hamiltonian system
(with $\eps$ as the perturbation parameter)
in the vicinity of a \emph{simple resonance},
after carrying out one step of resonant normal form,
and a~rescaling that gives rise to the fast frequencies~(\ref{eq:omega_eps})
(see \cite{Niederman00,DelshamsG01} for details).
The parameters $\eps$ and $\mu$ should not be regarded as independent,
but linked by a relation of the type $\mu=\eps^p$.

Notice that the unperturbed system $H_0$
consists of the pendulum given by $P(x,y)= y^2/2 + \cos x -1$ and
$2$~rotors with fast frequencies:
$\dot{\varphi}= \omega_\varepsilon + \Lambda I$, $\dot{I}=0$.
The pendulum has a hyperbolic equilibrium at the origin,
and the (upper) separatrix can be parametrized by
$(x_0(s),y_0(s))=(4 \arctan\ee^s,2/\cosh s)$, $s\in\R$.
The rotors system $(\varphi, I)$ has the solutions
$\varphi = \varphi_0+(\omega_\varepsilon + \Lambda I_0)\,t$, $I= I_0$.
Consequently, $H_0$ has a $2$-parameter family of $2$-dimensional
whiskered invariant tori, with coincident stable and unstable whiskers.
Among the family of whiskered tori,
we will focus our attention on the torus located at $I=0$,
whose frequency vector is $\omega_\varepsilon$ as in~(\ref{eq:omega_eps}),
in our case a frequency vector of constant type.
We also assume the condition of \emph{isoenergetic nondegeneracy}
\beq
    \det \left(
    \begin{array}{cc}
    \Lambda & \omega\\
    \omega^\top & 0
    \end{array}
    \right) \neq 0.
\label{eq:isoenerg}
\eeq
When adding the perturbation $\mu H_1$,
the \emph{hyperbolic KAM theorem} can be applied
(see for instance \cite{Niederman00})
thanks to the Diophantine condition~(\ref{eq:DiophCond})
and the isoenergetic nondegeneracy~(\ref{eq:isoenerg}).
For $\mu$ small enough, the whiskered torus persists
with some shift and deformation, as well as its local whiskers.

In general, for $\mu\ne0$ the (global) whiskers
do not coincide anymore, and one can introduce a \emph{splitting function}
giving the distance between the whiskers
in the directions of the action coordinates $I\in\R^2$:
denoting $\J^\pm(\theta)$ parameterizations
of a transverse section of both whiskers,
one can define \ $\M(\theta):=\J^-(\theta)-\J^+(\theta)$, \ $\theta\in\T^2$
\ (see \cite[\S5.2]{DelshamsG00}, and also~\cite{Eliasson94}).
This function turns out to be the gradient
of the (scalar) \emph{splitting potential}:
\ $\M(\theta)=\nabla\Lc(\theta)$.
\ Notice that the nondegenerate critical points of~$\Lc$ correspond
to simple zeros of $\M$ and give rise
to transverse homoclinic orbits to the whiskered torus.

Applying the Poincar\'e--Melnikov method,
the first order approximation~(\ref{eq:Melniapprox})
is given by the (vector) \emph{Melnikov function} $M(\theta)$,
which is the gradient of the \emph{Melnikov potential}:
\ $M(\theta)=\nabla L(\theta)$.
\ The latter one can be defined by integrating
the perturbation $H_1$ along
a trajectory of the unperturbed homoclinic manifold,
starting at the point of the section $s=0$ with a given phase~$\theta$:
\beq\label{eq:L}
  L(\theta) = - \int_{-\infty}^{\infty}
  [h(x_0(t))-h(0)] f(\theta +\omega_\varepsilon t)\,dt.
\eeq

In order to emphasize the r\^ole played by the arithmetic properties
of the frequencies, we have chosen for the perturbation the special form
given in~(\ref{eq:ham2}--\ref{eq:ham3}).
This form was considered in \cite{LochakMS03,DelshamsG04},
and allows us to deal with the Melnikov function and obtain
asymptotic estimates for the splitting.
Notice that the constant $\rho>0$ in the Fourier expansion
of $f(\varphi)$ in~(\ref{eq:ham3})
gives the complex width of analyticity of this function.
The phases $\sigma_k$ can be chosen arbitrarily for our purpose in this paper,
since we are focused in lower bounds for the splitting
(instead, some restriction on the phases has to be imposed in order to
the study the transversality of the splitting, as in \cite{DelshamsG04}).

Now we can formulate our main result, providing a lower estimate
for the \emph{maximal splitting distance}
in the case of frequencies of constant type.
We use the notation $|f|\succeq|g|$ if we can bound $|f|\ge a|g|$ with some
positive constant $a$ not depending on $\eps$, $\mu$.

\begin{theorem}[main result]\label{thm:main}
For the Hamiltonian system introduced
in~(\ref{eq:Hpert}), (\ref{eq:ham1}--\ref{eq:ham3})
with $3$ degrees of freedom,
satisfying the isoenergetic condition~(\ref{eq:isoenerg}),
assume that $\varepsilon$ is small enough and $\mu=\varepsilon^p$ with $p>3$.
If $\Omega$ in (\ref{eq:omega_eps}) is a number of constant type,
then the following lower bound holds:
\begin{equation}
 \max\limits_{\theta\in \mathbb{T}^2} |\M(\theta)|
\succeq \frac{\mu}{\varepsilon^{1/2}}\,\exp\pp{-\frac{C}{\eps^{1/4}}},
   \label{eq:lower}
\end{equation}
where $C=C(\Omega,\rho)$ is a positive constant,
defined in~(\ref{eq:defC}).
\end{theorem}

In the proof of this result, we put emphasis
on the constructive part of the proofs,
using the arithmetic properties of the frequencies
in order to provide a methodology which can be applied
to the case of numbers of constant type,
stressing the similarities and differences with the quadratic frequencies
considered in \cite{DelshamsGG13}.
We show that, for a given $\eps$ small enough,
the dominant harmonic of the Melnikov function $M(\theta)$
can be related to a resonant convergent,
and, consequently, we obtain an estimate
for the maximal value of this function.
In a further step, the first order approximation can be validated
showing that the dominant harmonics of the splitting
function $\M(\theta)$ correspond to the dominant harmonics
of the Melnikov function, as done in \cite{DelshamsG04}.

The paper is organized as follows. In Section~\ref{sect:arithm_nct}
we study the arithmetic properties of irrational numbers of constant type,
and in Section~\ref{sect:lower_bnd} we find, for frequencies
of constant type, the dominant harmonic of the splitting potential,
whose size allows us to provide
a lower bound for the maximal splitting distance,
as established in Theorem~\ref{thm:main}.

\section{Arithmetic properties of numbers of constant type}
  \label{sect:arithm_nct}

\subsection{Continued fractions and principal convergents}

Let $0<\Omega<1$ be an irrational number.
It is well-known that it has an infinite continued fraction
\beq\label{eq:cont_frac}
  \Omega=[a_1,a_2,a_3,\ldots]
  =\frac{1}{a_1+\dfrac{1}{a_2+\dfrac{1}{a_3+\cdots}}}\,,
  \qquad
  a_n\in\Z^+,\ n\ge1
  \quad \mbox{(and $a_0=0$)}.
\eeq
Its entries, which are integers $a_n\ge1$,
are called the \emph{partial quotients} of the continued fraction.
It is also well-known that the rational numbers
\ $\dfrac{p_n}{q_n}=[a_1,\ldots,a_n]$,\ $n\ge1$,
\ called the \emph{(principal) convergents} of $\Omega$,
provide successive best rational approximations to~$\Omega$.
Thus, if we consider the vectors $w(n):=(q_n,p_n)$, we obtain approximations to
the direction of the vector $\omega=(1,\Omega)$
(see, for instance, \cite{Schmidt80} and \cite{Lang95}
as general references on continued fractions).

Now, in order to obtain approximations
to the orthogonal line~$\langle\omega\rangle^\bot$,
i.e.~to the \emph{quasi-resonances} of $\omega$, we introduce
the following sequence of vectors, that we call
the \emph{``resonant convergents''}:
\begin{equation}\label{eq:vl}
  v(n):=(-p_n,q_n).
\end{equation}

We see from the standard recurrences
\[
  \begin{array}{lll}
    q_n=a_n q_{n-1}+q_{n-2}\,,   &\qquad q_0=1,       &\quad q_{-1}=0,
  \\[4pt]
    p_n=a_n p_{n-1} + p_{n-2}\,, &\qquad p_0=a_0=0\,, &\quad p_{-1}=1,
  \end{array}
\]
that the vectors $w(n)$ and $v(n)$ are also given by
recurrence relations:
\bea
  \label{eq:wlvl1}
  &&w(n)=a_n w(n-1)+w(n-2),
  \qquad
  w(0)=(1,0), \ w(-1)=(0,1),
\\[4pt]
  \label{eq:wlvl2}
  &&v(n)=a_n v(n-1)+v(n-2),
  \qquad
  \ \ v(0)=(0,1), \ \ v(-1)=(-1,0).
\eea
The following result provides alternative expressions for $w(n)$ and $v(n)$
in terms of products of unimodular matrices.
We point out that similar products for $w(n)$
appear in \cite{LiardetS98,Stambul00}.

\begin{proposition}\label{prop:wlvl}
Let $0<\Omega<1$ be irrational.
We define $A_m=\symmatrix{a_m}10$, $m\geq 1$, where $a_m$ are the partial
quotients of the continued fraction~(\ref{eq:cont_frac}).
Then, for $n\geq 1$ one has:
\begin{itemize}
\item[\rm(a)] $w(n) =  A_1\cdots A_n w(0)$, \quad $w(0)=(1,0)$,
\vspace{4pt}
\item[\rm(b)] $v(n) = (-1)^n A_1^{-1}\cdots A_n^{-1} v(0)$, \quad $v(0)=(0,1)$.
\end{itemize}
\end{proposition}

\proof
It is enough to prove the relations
\beq\label{eq:prodA}
  A_1 \cdots A_n = \mmatrix{q_n}{q_{n-1}}{p_n}{p_{n-1}},
  \qquad
  A_1^{-1}\cdots A_n^{-1} = (-1)^n \mmatrix{p_{n-1}}{-p_n}{-q_{n-1}}{q_n},
\eeq
since the right multiplication of each equation with $w(0)$ and $v(0)$,
respectively, implies the assertions~(a) and~(b).
To prove the equalities~(\ref{eq:prodA}),
we use induction and
use the recurrence relations~(\ref{eq:wlvl1}--\ref{eq:wlvl2}).
For $n=1$, the equalities are easily verified.
Assuming that the equalities are true for~$n$
we prove them for $n+1$:
\bean
  &&A_1 \cdots A_n A_{n+1}
  =\mmatrix{q_n}{q_{n-1}}{p_n}{p_{n-1}} \symmatrix{a_{n+1}}10
  =\mmatrix{q_{n+1}}{q_n}{p_{n+1}}{p_n},
\\[6pt]
  &&A_1^{-1}\cdots A_n^{-1} A_{n+1}^{-1}
  =(-1)^n\mmatrix{p_{n-1}}{-p_n}{-q_{n-1}}{q_n} \mmatrix{0}{1}{1}{-a_{n+1}}
  =(-1)^{n+1}\mmatrix{p_n}{-p_{n+1}}{-q_n}{q_{n+1}}.
\eean
\qed

\subsection{Numbers of constant type}

In this paper we are interested in the study of arithmetic properties of
the \emph{numbers of constant type}, i.e.~irrational numbers whose
continued fraction has bounded partial quotients.
Such numbers are also called \emph{badly approximable}.
We refer to \cite{Schmidt80}, \cite{Lang95} and to the survey \cite{Shallit92}
as general references.

Some explicit examples of non-quadratic numbers of constant
type appear in \cite{Shallit79,Shallit92}.
For instance, for the number
\beq\label{eq:shallit}
  \Omega=2\sum_{k=1}^\infty 2^{-2^k}
  =[1,1,1,2,1,1,1,1,1,1,1,2,1,1,1,\ldots]
  \simeq0.632843018043786,
\eeq
all partial quotients are~1 or~2,
and can be obtained using a recurrence relation,
given in \cite[Th.~11]{Shallit79}.

The property of being a number of constant type is equivalent
to the Diophantine condition with the minimal exponent~(\ref{eq:DiophCond}).
In view of this, we define as in \cite{DelshamsG03}
the following \emph{``numerators''}:
\beq\label{eq:numerators}
  \gamma_k := |\langle k, \omega\rangle|\cdot|k|
  \qquad
  \forall k\in\Z^2\setminus\pp0
\eeq
(for integer vectors, we use the norm $\abs\cdot=\abs\cdot_1$,
i.e.~the sum of absolute values of the components of the vector).
Our aim is to study the integer vectors $k$ which give
the smallest values $\gamma_k$.
We denote
\begin{equation}
\gamma^*:=\liminf_{|k|\to \infty} \gamma_k > 0.
\label{eq:gamast}
\end{equation}
We can also define analogous numerators
from the expression~(\ref{eq:DiophcondC}) of the Diophantine condition,
as well as their associated asymptotic value
(see also \cite{Cassels57, Schmidt80}):
\begin{equation}\label{eq:c_ast}
  \nu_q:=q\norm{q\Omega} \quad \forall q\geq 1,
  \qquad
  \nu^*:=\ds\liminf_{q\to\infty}\nu_q > 0,
\end{equation}
where we denote $\norm a:=\abs{a-\rint(a)}=\ds\min_{p\in\Z}\abs{a-p}$,
i.e.~the distance to the closest integer.
It is easy to check that $\nu^*=\gamma^*/(1+\Omega)$.
Indeed, writing $k=(-p,q)$ in~(\ref{eq:DiophCond}), we have that,
for a fixed $q\geq 1$, the small divisor
$|\langle k, \omega \rangle| = |q \Omega -p|$ is minimized
for $p=\rint(q\Omega)$, and we have $|k|=p+q \approx (1+\Omega) q$,
as $q\to \infty$.

\paragr{Quadratic numbers}
We give a brief summary of the results of \cite{DelshamsG03},
concerning quadratic irrational numbers
(i.e.~real roots of quadratic polynomials with integer coefficients).
Notice that all quadratic numbers belong
to the class of numbers of constant type.
Indeed, this is a consequence of the well-known fact
that for a quadratic number $\Omega$,
the continued fraction~(\ref{eq:cont_frac}) is eventually periodic,
with some period $m\geq1$,
that is, there exists $l\ge1$ such $a_{n+m}=a_n$ for $n\geq l$.

A technique for studying resonances for quadratic
frequencies was developed in \cite{DelshamsG03}, where
the periodicity of the continued fraction of quadratic numbers
was used to construct a  unimodular matrix $T$ having the vector
$\omega=(1,\Omega)$ as an eigenvector with eigenvalue $\lambda>1$.
Clearly, the matrix $T$ provides approximations to the direction of $\omega$.
Then, the associated quasi-resonances are given
by the matrix $U:=(T^{-1})^\top$.
In fact, the study can be restricted to the case of purely periodic numbers
$\Omega=[\ol{a_1,\ldots,a_m}]$.
Otherwise for a quadratic number with non-purely periodic continued fraction
$\hat\Omega=[b_1,\ldots,b_l, \ol{a_{1},\ldots,a_{m}}]=[b_1,\ldots,b_l,\Omega]$
a linear change given by a unimodular matrix can be done between
$\omega=(1,\Omega)$ and $\hat\omega=(1,\hat\Omega)$.
If $\Omega\in(0,1)$ is a quadratic irrational number
with a purely $m$-periodic continued fraction,
then  $T=A_1\cdots A_m$ and $U=A_1^{-1}\cdots A_m^{-1}$,
where the matrices $A_n$ have been
introduced in Proposition~\ref{prop:wlvl}.
All the integer vectors $k\in\Z^2\setminus\pp0$
with $\abs{\scprod k\omega}<1/2$
can be subdivided into \emph{resonant sequences}:
\begin{equation}
s(j,n) := U^n k^0(j),
\qquad n=0,1,2,\ldots
\label{eq:sjn}
\end{equation}
where the initial vector $k^0(j)=(-\rint(j\Omega),j)$, $j\in \Z^+$, satisfies
\begin{equation}\label{eq:primit}
 \frac{1}{2\lambda} < |\langle k^0(j), \omega\rangle| < \frac{1}{2}.
\end{equation}
For each $j\in \Z^+$ satisfying~(\ref{eq:primit}),
it was proved in \cite[Th.~2]{DelshamsG03} (see also \cite{DelshamsGG13})
that, asymptotically, the resonant sequence $s(j,n)$
exhibits an geometric growth as $n\to\infty$, with ratio $\lambda$, and that
the sequence of the numerators $\gamma_{s(j,n)}$ has a limit $\gamma^*_j$.
Since the lower bounds for $\gamma^*_j$, also provided in \cite{DelshamsG03},
are increasing in $j$, we can select the minimal of them
corresponding to some $j_0$ and, thus,
get the value for $\gamma^*$ defined in~(\ref{eq:gamast}).
The integer vectors of the corresponding sequence $s(j_0, n)$
are called \emph{the primary resonances},
and the \emph{secondary resonances} are  the integer vectors
belonging to any of the remaining sequences $s(j,n)$, $j\neq j_0$.

It is also easy to deduce that $U v(n)=(-1)^m v(n+m)$,
where $v(n)$ is defined in (\ref{eq:vl}).
This implies that for quadratic frequencies
the sequence of resonant convergents $v(n)$
is divided into $m$~resonant sequences~(\ref{eq:sjn})
(the sign is not relevant).

\paragr{Non-quadratic numbers of constant type}
If $\Omega$ is a number of constant type, non-quadratic,
there is no periodicity in its continued fraction and, hence, we cannot
construct a matrix like $U$, or a sequence analogous to the primary resonances.
Alternatively, we can use in this case the sequence of
resonant convergents $v(n)$, defined in (\ref{eq:vl}),
to obtain some similar results.

In the next result we provide an upper and
a lower bound for $\nu^*$ in~(\ref{eq:c_ast}), and, hence, for $\gamma^*$
in~(\ref{eq:gamast}), and we show that it can be obtained by restricting~$q$
to the denominators $q_n$ of the convergents.

\begin{lemma}\label{lm:c_ast}
For the numbers $0<\Omega<1$ of constant type,
we have
$$
\frac{1}{3} \leq \nu^*=\ds\liminf_{n\to\infty} \nu_{q_n}
\leq \frac{1}{\sqrt{5}}\,,
\qquad
\ds\limsup_{n\to\infty} \nu_{q_n}\leq 1.
$$
\end{lemma}

\proof
We use several results in \cite[\S I.5]{Schmidt80}
(namely, Theorems~I.5A, I.5B, I.5C and Lemma~I.3E),
concerning the properties of the convergents
of any irrational number. On one hand,
for a given $q$ and $p=\rint(q\Omega)$,
if $\abs{q\Omega-p}<1/2q$ and $p/q$ is a reduced fraction,
i.e.~$\nu_q<1/2$, then $q$ is a convergent $q_n$\,, 
and an infinite number of convergents satisfies such inequality 
(otherwise, if it is not a reduced fraction then we can write $p/q=p_n/q_n$
for some convergent, and $\nu_{q_n}<\nu_q$).
This implies that the limit in~(\ref{eq:c_ast}) can be restricted
to the convergents $q_n$\,.
In fact, an infinite number of convergents satisfy the sharper inequality
$\norm{q_n\Omega}<1/(\sqrt{5}\,q_n)$, 
which gives the upper bound $\nu^*\le1/\sqrt5$\,.
On the other hand, since all convergents satisfy
the inequality $\norm{q_n\Omega}<1/q_n$\,, 
we get the upper bound \ $\ds\limsup_{n\to\infty} \nu_{q_n}\le1$.
\ Finally, a classical Markoff's theorem \cite{Markoff79}
implies that $\nu^*\geq1/3$.
\qed

Now we define
\begin{equation}\label{eq:nct_E}
  E=E(\Omega)
  :=\p{\frac1{\nu^*}\limsup_{n\to\infty}\nu_{q_n}}^{1/2}
  =\p{\frac{\ds\limsup_{n\to\infty}\nu_{q_n}}{\ds\liminf_{n\to\infty}\nu_{q_n}}}^{1/2},
\end{equation}
which is a finite number (the bounds of Lemma~\ref{lm:c_ast}
imply that $1\leq E\leq\sqrt3$),
and the bounds
$\nu^*\leq \nu_{q_n} \leq E^2 \nu^*$ hold asymptotically, as $n\to \infty$.
(We use an expression like ``\,$a_n\le b_n$ as $n\to\infty$\,''
if $\ds\limsup_{n\to\infty}(a_n/b_n)\le1$.)

Consequently, the smallest numerators $\gamma_k$ can be found among
the resonant convergents $k=v(n)$ introduced in~(\ref{eq:vl}).
Moreover, we have the following asymptotic bounds, as $n\to\infty$:
\beq\label{eq:boundsgamma}
  \gamma^*\leq \gamma_{v(n)} \leq E^2 \gamma^*.
\eeq

The following result provides a \emph{geometric lower and upper bound}
for the convergents $q_n$\,,
generalizing in some sense the geometric asymptotic estimate
of resonant sequences of a quadratic number
given in \cite[Th.~2]{DelshamsG03}.

\begin{proposition}\label{prop:geom}
Let $\Omega=[a_1,a_2,\ldots]$ be a number of constant type,
and define \ $M=M(\Omega):=1+\ds\max_na_n$\,. For any $n\ge2$, one has:
\[
  1+\frac1M\le\frac{q_n}{q_{n-1}}\le M.
\]
\end{proposition}

\proof
This is a simple consequence of a well-known general formula:
\ $\dfrac{q_n}{q_{n-1}}=a_n+[a_{n-1},\ldots,a_1]$
\ (see \cite[Lemma~\S I.3F]{Schmidt80}),
which implies the inequalities
\ $\ds a_n+\frac1{a_{n-1}+1}=a_n+[a_{n-1},1]\le\frac{q_n}{q_{n-1}}
\le a_n+[a_{n-1},\infty]\le a_n+1$.
\qed

\section{Lower bounds for the maximal splitting distance}\label{sect:lower_bnd}

This section is devoted to the proof of Theorem~\ref{thm:main},
giving a lower bound for the maximal splitting distance
between the invariant manifolds of the whiskered torus,
when the frequency ratio is a number of constant type.
We start with a brief description of our approach.
As said in the introduction, a measure for the splitting distance
on a transverse section is given by the splitting function, which is the
gradient of the splitting potential:
$\M(\theta)=\nabla\Lc(\theta)$, $\theta\in\T^2$.
In order to obtain a lower bound for the splitting function,
we consider the Fourier expansion of the splitting potential:
\[
  \Lc(\theta)=\hardsum{k\in \mathbb{Z}^2\setminus\{0\}}{k_2\ge0}
  \Lc_k\cos(\langle k,\theta\rangle - \tau_k).
\]
For the splitting function, its vector Fourier coefficients $\M_k$
are related to the scalar coefficients of the splitting potential:
$\abs{\M_k}=\abs k\Lc_k$\,.

The Poincar\'e--Melnikov method provides the first order
approximation~(\ref{eq:Melniapprox}) in terms of the Melnikov function,
which is the gradient of the splitting potential:
$M(\theta)=\nabla L(\theta)$.
The Fourier coefficients $L_k$ of the Melnikov potential can be
computed explicitly, which allows us to detect the dominant harmonic,
corresponding to some $k=S(\eps)$ for any sufficiently small $\eps$.
Such dominance, which holds for the Poincar\'e--Melnikov
approximation and is found in a constructive way,
can also be validated for the whole splitting function
when the error term in~(\ref{eq:Melniapprox}) is added,
in our singular case $\mu=\varepsilon^p$,
analogously to \cite{DelshamsG04}
(where the case of the golden number was considered).

In this way, we obtain an exponentially small asymptotic estimate
for the coefficient dominant harmonic,
as well as an exponentially small upper bound
for the sum of all other harmonics, showing that the dominant harmonic
is large enough to ensure that the maximal splitting distance
can be approximated by the size of this coefficient:
\beq\label{eq:dominantcoef}
  \max_{\theta\in\T^2}\abs{\M(\theta)}\approx\abs{\M_S}
  \quad\textrm{as $\eps\to0$},
  \qquad
  S=S(\eps).
\eeq
(We use the expression ``\,$f(\eps)\approx g(\eps)$ as $\eps\to0$\,''
if $\ds\lim_{\eps\to0}(f(\eps)/g(\eps))=1$.)
The asymptotic estimate for $\abs{\M_S}$ is given
in Proposition~\ref{prop:Mdom}(a),
in terms of an oscillating function $h_1(\eps)$.
If $\Omega$ is a quadratic irrational number, the function $h_1(\eps)$
is periodic with respect to $\ln\eps$ (see \cite{DelshamsG03,DelshamsGG13}).
Instead, in the case of a number of constant type $\Omega$
under consideration, the function $h_1(\eps)$ has a more complicated behavior
but can be bounded from below and above, which gives rise
to exponentially small upper bounds
(as in \cite{Simo94}, \cite{DelshamsGS04} and other works)
and to exponentially small lower bounds
as in the statement of Theorem~\ref{thm:main}.

\subsection{Dominant harmonics of the splitting potential}\label{sect:dominant}

We put our functions $f$ and $h$ defined in~(\ref{eq:ham3})
into the integral (\ref{eq:L}) and, calculating it by residues,
we get the Fourier expansion of the Melnikov potential:
\[
  L(\theta)=\hardsum{k\in \mathbb{Z}^2\setminus\{0\}}{k_2\ge0}
  L_k \cos(\langle k, \theta\rangle -\sigma_k),
  \qquad
  L_k = \frac{2\pi |\langle k, \omega_\varepsilon\rangle|
  \,\ee^{-\rho |k|}}{\sinh |\frac{\pi}{2}
  \langle k, \omega_\varepsilon\rangle|}\,.
\]
Using~(\ref{eq:omega_eps}) and~(\ref{eq:numerators}),
we can present the coefficients in the form
\begin{equation}
\label{eq:alphabeta}
L_k = \alpha_k\,\ee^{- \beta_k},
\qquad
\alpha_k \approx \frac{4 \pi\gamma_k}{|k|\sqrt{\varepsilon}}\,,
\quad
\beta_k =\rho |k| + \frac{\pi \gamma_k}{2 |k|\sqrt{\varepsilon}}\,,
\end{equation}
where an exponentially small term has been neglected in the denominator
of $\alpha_k$.
For any given $\eps$, the harmonics with largest coefficients $L_k(\eps)$
correspond essentially to the smallest exponents
$\beta_k(\eps)$. Thus, we have to study the dependence on $\eps$ of
such exponents.

With this aim, we introduce for any $X$, $Y$
the function
\[
  G(\eps;X,Y):=
  \frac{Y^{1/2}}2
  \pq{\p{\frac\eps X}^{1/4}+\p{\frac X\eps}^{1/4}},
\]
having its minimum at $\varepsilon=X$,
with the minimum value $G(X;X,Y)=Y^{1/2}$.
Then, the exponents $\beta_k(\eps)$ in~(\ref{eq:alphabeta}) can be
presented in the form
\beq\label{eq:gk_quad}
  \beta_k(\eps) = \frac{C_0}{\eps^{1/4}}\,g_k (\eps),
  \qquad
  g_k (\eps):=G(\eps;\eps_k,\tl\gamma_k),
  \qquad
  C_0:=(2\pi\rho\gamma^*)^{1/2},
\eeq
where
\beq\label{eq:C0}
  \eps_k:=D_0\,\frac{\tl\gamma_k^{\,2}}{\abs k^4}\,,
  \qquad
  \tl\gamma_k := \frac{\gamma_k}{\gamma^*}\,,
  \qquad
  D_0:=\p{\frac{\pi\gamma^*}{2\rho}}^2.
\eeq
Notice that, to define the \emph{``normalized numerators''} $\tl\gamma_k$\,,
we have taken into account the limit numerator~$\gamma^*$
introduced in~(\ref{eq:gamast}).
Consequently, for all $k$ we have
\ $\beta_k(\eps)\geq \dfrac{C_0\tl\gamma_k^{1/2}}{\varepsilon^{1/4}}$\,,
\ which provides an asymptotic estimate for
the maximum value of the coefficient $L_k(\eps)$ of each harmonic.
We point out that similar estimates of the size of a given harmonic
from the arithmetic properties of frequencies, for the case of constant type,
were already obtained in \cite{RudnevW98}.

We define, for any given $\eps$, the functions $h_1(\eps)$ and $h_2(\eps)$
as the first and the second minima of the values
$g_k(\eps)$, $k\in\Z^2\setminus\pp0$,
and we denote $S=S(\eps)$ the integer vector that gives the first minimum:
\beq\label{eq:h12}
  h_1(\eps):=\min_kg_k(\eps)=g_S(\eps),
  \qquad
  h_2(\eps):=\min_{k\ne S}g_k(\eps).
\eeq
The function $h_1(\eps)$ is continuous, and provides an estimate of the size
of the most dominant coefficient of the Melnikov potential.
The function $h_2(\eps)$ is also continuous,
and we have $h_2>h_1$ excepting at some isolated values of $\eps$
where $h_1$ and $h_2$ coincide
because of a change in the vector $S(\eps)$ giving the dominant harmonic.
Notice also that $S(\eps)$ remains constant in each interval between
two of such consecutive values of $\eps$.

\begin{figure}[!b]
  \centering
  \includegraphics[width=0.83\textwidth]{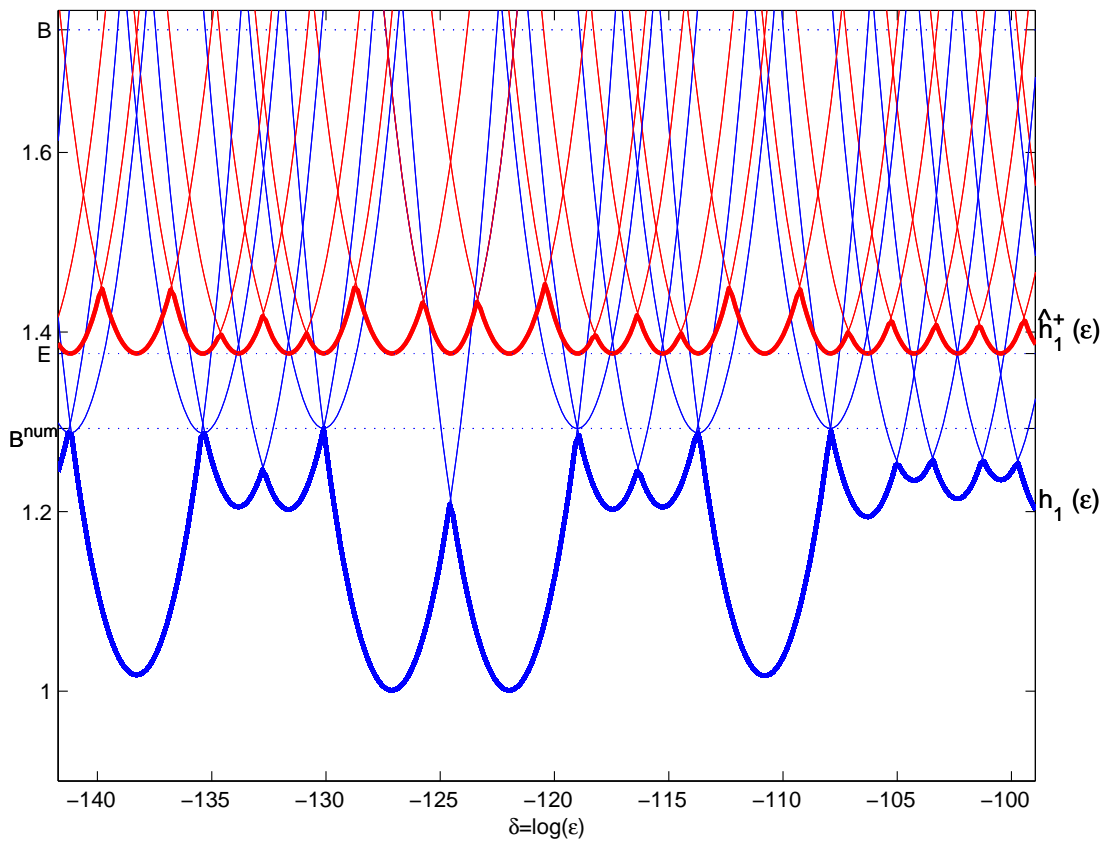}
  \caption{\small\emph{Graphs of the functions
    $\wh g_n(\varepsilon)$, $h_1(\varepsilon)$ (in thin and thick blue) and
    $\wh g_n^+(\varepsilon)$, $\wh h_1^+(\varepsilon)$ (in thin and thick red)
    using a logarithmic scale for $\eps$,
    for the frequency ratio $\Omega$ given in~(\ref{eq:shallit}).
    (Constants: $M=3$, $E\approx1.3761$, $B\approx1.7366$,
    $B^{\emph{num}}\approx1.2925$.)}}
\label{fig:nct2}
\end{figure}

In the particular case of quadratic frequencies,
the function $h_1(\eps)$ is periodic in $\ln\eps$
(see \cite{DelshamsG03,DelshamsGG13}),
due to a periodicity in the graphs of the functions $g_k(\eps)$,
which can be deduced from the asymptotic geometric growth of
the resonant sequences~(\ref{eq:sjn}).
Unfortunately, in the more general case of non-quadratic
frequencies of constant type,
the graphs of $g_k(\eps)$ do no exhibit a periodicity in $\ln\eps$.
In general, it is hard to provide an analytic description
of the function $h_1(\eps)$, but we can obtain an upper bound for it,
as $\eps\to0$ (see Section~\ref{sect:h1_bnd}),
which gives rise to a lower bound for the dominant harmonic
and, in view of~(\ref{eq:dominantcoef}),
also for the maximal splitting distance.
As an illustration, the function $h_1(\eps)$
is represented in Figure~\ref{fig:nct2}
for the concrete case of the frequency ratio $\Omega$
introduced in~(\ref{eq:shallit}).

Now we give, from the function $h_1(\eps)$,
an estimate for the dominant harmonic in the Fourier expansion
of the splitting function, as well as for the sum of all other harmonics.
We use the notation $f\sim g$ if $c_1|g|\le|f|\le c_2|g|$ with some positive
constants $c_1$, $c_2$ not depending on~$\eps$, $\mu$.

\begin{proposition}\label{prop:Mdom}
For $\eps$ small enough and $\mu=\varepsilon^p$ with $p>3$, one has:\\
\btm
\item[\rm(a)]
$\ds\abs{\M_S}
 \sim\mu\abs{S}L_S
 \sim\frac{\mu}{\varepsilon^{1/2}}
   \exp\pp{-\frac{C_0h_1(\eps)}{\eps^{1/4}}}$;
\vspace{4pt}
\item[\rm(b)]
$\ds\sum_{k\ne S}\abs{\M_k}
 \sim\frac{\mu}{\varepsilon^{1/2}}
   \exp\pp{-\frac{C_0h_2(\eps)}{\eps^{1/4}}}$.
\etm
\end{proposition}

\proof
The proof is similar to the analogous result in \cite[Prop.~4]{DelshamsGG13}.
By the Poincar\'e--Melnikov method~(\ref{eq:Melniapprox}),
the coefficients of the splitting function can be approximated by
the expression
\[
  \abs{\M_k}=\abs k\Lc_k\sim\mu\abs kL_k=\mu\abs k\alpha_k\,\ee^{-\beta_k},
\]
where we have neglected the error term of~(\ref{eq:Melniapprox})
in this first approximation,
and we have used the expression~(\ref{eq:alphabeta}) of the coefficients of
the Melnikov potential.
As mentioned previously,
the main behavior of the coefficients~$L_k(\eps)$ is given by
the exponents $\beta_k(\eps)$, which have been written in~(\ref{eq:gk_quad})
in terms of the functions~$g_k(\eps)$.
In particular, the coefficient $L_S$\,,
associated to the dominant harmonic $S=S(\eps)$,
can be expressed in terms of the function $h_1(\eps)$
introduced in~(\ref{eq:h12}).
Thus, the exponential factor in~(a) is directly given by~$\ee^{-\beta_S}$.
On the other hand, the polynomial factor comes from
the estimate \ $\abs S\alpha_{S}\sim1/\sqrt\eps$, \ which comes directly from
the approximation of $\alpha_S$ given by~(\ref{eq:alphabeta}),
using also the estimate $\gamma_S\sim1$
(this is a consequence of Lemma~\ref{lm:h1bounds},
see the first remark in the next section).
In this way, we have obtained an asymptotic estimate for the
size $\abs{M_S}$ of the dominant coefficient of the Melnikov function.
To complete the proof of part~(a), one has to see
that the same estimate is valid for $\M_S$,
i.e.~when the error term in the Poincar\'e--Melnikov
approximation ~(\ref{eq:Melniapprox}) is not neglected,
in our singular case $\mu=\eps^p$.
We omit the details for this step, since it can be worked out
as it was done in \cite[Lemma~5]{DelshamsG04} for the case
of the golden number ($\Omega=(\sqrt5-1)/2$),
using upper bounds for the error term provided in \cite[Th.~10]{DelshamsGS04}.

The proof of part~(b) is carried out in similar terms.
For the second dominant harmonic, we get an exponentially small estimate
with the function $h_2(\eps)$, defined in~(\ref{eq:h12}).
This estimate is also valid if one considers the whole sum in~(b),
since for any given $\eps$ the terms of this sum
can be bounded by geometric series
and, hence, it can be estimated by its dominant term
(see \cite[Lemma~4]{DelshamsG04} for more details).
\qed

\subsection{Upper and lower bounds for $h_1(\eps)$}\label{sect:h1_bnd}

To conclude the proof of Theorem~\ref{thm:main}, we provide bounds for the
function $h_1(\eps)$. In particular, an upper bound for this function
gives rise to a lower bound for the splitting function,
in view of the approximation~(\ref{eq:dominantcoef})
and Proposition~\ref{prop:Mdom}(a).

We are going to show that, to provide an upper bound for $h_1(\eps)$,
we can restrict our study of the minimum in~(\ref{eq:h12})
to the harmonics associated to resonant convergents $v(n)$,
defined in~(\ref{eq:vl}),
which play an analogous r\^ole as the primary resonances
for quadratic frequencies considered in \cite{DelshamsG03}.
We are going to use the results on arithmetic properties of
frequencies of constant type (see Section~\ref{sect:arithm_nct}).
Thus, we define
\[
  \wh h_1(\eps):=\min_n\wh g_n(\eps).
  \qquad
  \wh g_n(\eps):=g_{v(n)}(\eps)=G(\eps;\eps_{v(n)},\tl\gamma_{v(n)}),
\]
We know from~(\ref{eq:boundsgamma}) that
the minimum values of the functions $\wh g_n(\eps)$ satisfy the bounds:
$1\le\tl\gamma_{v(n)}^{1/2}\le E$ (asymptotically, as $n\to\infty$),
where the constant $E=E(\Omega)$ has been defined in~(\ref{eq:nct_E}).
With this in mind, we also define the functions
\begin{equation}\label{eq:gnplus_nct}
  \wh h_1^+(\eps):=\min_n\wh g_n^+(\eps),
  \qquad
  \wh g_n^+(\eps):=G(\eps;\eps_{v(n)},E^2),
\end{equation}
which are represented in Figure~\ref{fig:nct2}
for the concrete case of the frequency ratio $\Omega$
inroduced in~(\ref{eq:shallit}).
It is clear that \ $h_1(\eps)\le\wh h_1(\eps)\le\wh h_1^+(\eps)$
\ for any~$\eps$.
Since all the functions $\wh g^+_n(\eps)$ have the same
minimum value~$E$, it is much simpler to study the behavior
of function $\wh h^+_1(\eps)$ and, in particular,
to provide an upper bound for it.
As shown in the next lemma, this upper bound is given by the constant
\begin{equation}\label{eq:A1}
 B=B(\Omega):=\dfrac E2\pq{(EM)^{1/2}+(EM)^{-1/2}},
\end{equation}
where $E$ and $M$ have been defined in~(\ref{eq:nct_E})
and in Proposition~\ref{prop:geom}, respectively.

\begin{lemma}\label{lm:h1bounds}
The function $h_1(\eps)$ defined in (\ref{eq:h12})
satisfies the asymptotic bounds:
\ $1 \leq h_1(\eps) \leq B$, \ as $\eps\to0$.
\end{lemma}

\proof
The lower bound is very simple: the function $h_1(\eps)$
is defined from the functions $g_k(\eps)$ in~(\ref{eq:gk_quad}),
whose minimum values are $\tl\gamma_k^{1/2}$.
Recall that the normalized numerators were introduced in~(\ref{eq:C0}),
and satisfy the asymptotic bound $\tl\gamma_k\ge1$, as $\abs k\to\infty$.

Now we are going to obtain an upper bound for the function $\wh h_1^+(\eps)$
defined in~(\ref{eq:gnplus_nct}), and hence for $h_1(\eps)$.
According to~(\ref{eq:gk_quad}) and~(\ref{eq:C0}),
each function $\wh g^+_n(\eps)$ takes its minimum value at
\ $\eps_{v(n)}=D_0\tl\gamma_{v(n)}^{\,2}/\abs{v(n)}^4$.
\ Using the asymptotic bounds~(\ref{eq:boundsgamma}), as well as
the geometric growth of the convergents described
in Proposition~\ref{prop:geom}, we see that
\ $\eps_{v(n-1)}/\eps_{v(n)}\le(EM)^4$, \ as $n\to\infty$.
In other words, using a logarithmic scale for $\eps$
(as in Figure~\ref{fig:nct2}), the distance between the
minimum points $\eps_{v(n-1)}$ and $\eps_{v(n)}$ of two consecutive functions,
$\wh g^+_{n-1}$ and $\wh g^+_n$\,, is $\le4\ln(EM)$.
Their value at the ``middle point''
(which provides an approximate intersecting point) is
the constant $B$ defined in~(\ref{eq:A1}),
which provides an upper bound, as $\eps\to0$, for $\wh h^+_1(\eps)$
and also for~$h_1(\eps)$.
\qed

\bremarks
\item
In the example considered in Figure~\ref{fig:nct2},
for all values of $\eps$ the dominant harmonic $S=S(\eps)$ is given
by a resonant convergent $v(N)$, $N=N(\eps)$,
and we have $h_1(\eps)=\wh h_1(\eps)\le\wh h_1^+(\eps)$.
In other examples, it may happen that for some small intervals of $\eps$
the dominant harmonic $S$ is a non-convergent vector,
having $h_1(\eps)<\wh h_1(\eps)$ in such intervals.
In any case, the minimum value $\tl\gamma_S^{1/2}$ of the function $g_S(\eps)$
associated to the dominant harmonic $S=S(\eps)$ is always contained
in the interval $[1,B(\Omega)]$, hence we can write
$\gamma_S=\tl\gamma_S\cdot\gamma^*\sim1$, an estimate that was used in the
proof of Proposition~\ref{prop:Mdom}.
\item
The upper bound $h_1(\eps)\le B$ bound is not sharp,
since it considers the worst possible case in the bounds
of~(\ref{eq:boundsgamma}) and Proposition~\ref{prop:geom}.
A much sharper upper bound $h_1(\eps)\le B^{\textrm{num}}$
can be obtained numerically (see Figure~\ref{fig:nct2}).
\eremarks

Now we can complete the proof of Theorem~\ref{thm:main}.
Indeed, applying the upper bound $h_1(\eps)\le B$ to
the exponent in Proposition~\ref{prop:Mdom}(a)
and using~(\ref{eq:dominantcoef}), we get the lower bound~(\ref{eq:lower}),
and the constant in the exponent is given by
\beq\label{eq:defC}
  C(\Omega,\rho)=C_0(\Omega,\rho)\cdot B(\Omega),
\eeq
where $C_0$ and $B$ have been defined in~(\ref{eq:C0}) and~(\ref{eq:A1})
respectively.

\small

\def\noopsort#1{}

\end{document}